\documentclass[11pt]{amsart}
\usepackage{amssymb}
\newcommand{\cc}{\mathbb{C}}
\newcommand{\com}{\complement}

\newtheorem{theorem}{Theorem}
\newtheorem{lemma}[theorem]{Lemma}
\newtheorem{corollary}[theorem]{Corollary}
\theoremstyle{definition}
\newtheorem{example}[theorem]{Example}

\title[Complements of rationally convex sets] 
{Complements of rationally convex sets}

\author{E. S. Zeron}
\address{Depto. Matem\'aticas, CINVESTAV, Apartado 
Postal 14-740, M\'exico D.F., 07000, M\'exico.}
\email{eszeron@math.cinvestav.mx}

\date{\today}
\thanks{Research supported by Cinvestav and Conacyt 
(Mexico), and Universit\'e de Montr\'eal (Canada)}

\begin{document}
\begin{abstract}
The main objective of this paper is to show that the complement 
of a rational convex set in $\cc^n$ is $(n-2)$-connected.
\end{abstract}

\maketitle
\section{Introduction}

A compact subset $K$ of $\cc^n$ is said to be \textit{rationally}
(respectively \textit{polynomially}) convex if for every point $y_0$ 
in the complement $\cc^n\setminus{K}$ there exists a non-constant 
holomorphic polynomial $P$ such that $P(y_0)=1$ and $1\notin{P}(K)$ 
(respectively $\|P\|_K<1$). Notice that each polynomially convex set 
is rationally convex.

Rationally and polynomially convex sets play an extremely important 
role in complex analysis and approximation theory; we refer the reader 
to Stolzenberg~\cite{St} and Alexander and Wermer~\cite{AW}. Nevertheless, 
it is usually a difficult problem to decide whether a given compact set 
is rationally or polynomially convex. Therefore, the results that provide 
topological obstructions to rational or polynomial convexity are of special 
interest to complex analysis. Recently, Forstneri\v{c}~\cite{Fo} proved via 
Morse theory that the complement of a polynomially convex set is simply 
connected. The main objective of this paper is to show that the complement 
of a rationally convex set is simply connected as well.

\begin{theorem}[\textbf{Main}]\label{main}
Let $K$ be a compact rationally convex set in $\cc^n$, the 
complement $\cc^n\setminus{K}$ is simply connected for $n\geq3$. 
\end{theorem}

We prove Theorem~\ref{main} in the second section of this paper. Notice
that the complement of any compact rationally convex set $K$ in $\cc^n$ 
is arcwise connected whenever $n\geq2$. Actually, the complement of $K$ 
has only one unbounded arcwise connected component because $K$ is compact; 
and $\cc^n\setminus{K}$ has no bounded arcwise connected components 
because it is \textit{exhausted} by unbounded algebraic hypersurfaces
$Q^{-1}(1)\subset\cc^n$, with $Q$ a non-constant holomorphic polynomial 
of several variables. Thus, we need not specify any \textit{base} point 
for calculating the fundamental group of $\cc^n\setminus{K}$ in
Theorem~\ref{main}.

Moreover, the proof of Theorem~\ref{main} is based in the fact that the
complement $\cc^n\setminus{K}$ may also be \textit{exhausted} by smooth
algebraic hypersurfaces $\mathcal{H}\subset\cc^n$ of complex dimension
$n-1$, and so the relative pair $(\mathcal{H},\mathcal{H}\setminus{B}_R)$ 
is arcwise and simply connected for $n\geq3$ and $B_R$ an open ball in 
$\cc^n$ with centre in the origin and radius $R>0$ large enough; see 
Lefschetz theorem in \cite[p.~552]{Ha}. Suppose from now on that 
$\mathcal{H}$ is an algebraic hypersurface in $\cc^n$, with $n\geq4$, 
and recall Theorem~5.2 of \cite[p.~45]{Mi}. We have that the intersection 
$\mathcal{H}\cap\mathcal{S}_\delta$ is $(n-3)$-connected for every sphere
$\mathcal{S}_\delta\subset\cc^n$ with centre in $y_0\in\mathcal{H}$ and 
radius $\delta>0$ small enough. Thus, we may consider the question whether 
the intersection $\mathcal{H}\cap\mathcal{S}$ is also simply connected 
for any arbitrary sphere $\mathcal{S}\subset\cc^n$. A positive answer to 
previous question would shine light on a proof to the following possible  
extension of Theorem~\ref{main}:  The complement $\mathcal{S}\setminus{K}$ 
is simply connected for each rationally convex set $K$ contained in a 
sphere $\mathcal{S}\subset\cc^n$ with $n\geq4$. Unfortunately, the 
intersection of a smooth algebraic hypersurface $\mathcal{H}\subset\cc^n$ 
and a sphere $\mathcal{S}\subset\cc^n$ is not necessarily simply 
connected, as the following example shows,

\begin{example}
Consider the polynomial $Q_2(z)=z_1^2+z_2^2$ defined on $\cc^n$ with
$n\geq3$. It is easy to calculate that the hypersurface $Q_2^{-1}(1)$ 
is smooth, and the intersection of $Q_2^{-1}(1)$ with the unitary sphere
$\{\sum|z_k|^2=1\}$ is equal to the one-dimensional circumference 
defined by the equations: $z_k=0$ for $k\geq3$,
$$\Re(z_1)^2+\Re(z_2)^2=1\quad\hbox{and}\quad\Im(z_1)=\Im(z_2)=0.$$
\end{example}

In any case, it will be quite interesting to deduce whether the complement
$\mathcal{S}\setminus{K}$ is simply connected for each rationally convex
set $K$ contained in a sphere $\mathcal{S}\subset\cc^n$ with $n\geq4$. On
the other hand, Theorem~\ref{main} does not hold for the $2$-dimensional
space $\cc^2$ as the following simply example shows.

\begin{example}  
Consider the standard torus $\mathcal{T}^2$ in $\cc^2$ defined by
$|x|=|y|=1$. It is easy to see that $\mathcal{T}^2$ is rationally convex.  
Let $\Upsilon_0$ be a path parametrised by $t\mapsto(2,t)\sin(t)$ in the
real plane $\mathbb{R}^2$, for $0\leq{t}\leq\pi$. The path $\Upsilon_0$
begins at the origin, runs one time around the point $(1,1)$ and comes
back to the origin again. The compact set $E$ composed of those points
$(x,y)\in\cc^2$ with absolute values $(|x|,|y|)$ in $\Upsilon_0$ is 
then homeomorphic to $\Upsilon_0\times\mathcal{T}^2$ with the fibre
$\{0\}\times\mathcal{T}^2$ identified as a single point. Moreover, the
intersection of $\mathcal{T}^2$ and $E$ is empty, the set $E$ is a strong
deformation retract of the complement $\cc^2\setminus\mathcal{T}^2$, and
so, the first homotopy group of $\cc^2\setminus\mathcal{T}^2$ is equal 
to the integer numbers. \end{example}

We strongly recommend the works of Bredon~\cite{Be} and Spanier~\cite{Sp}
for references on homotopy theory. Besides, recalling the works of
Givental$'$~\cite{Gi}, Duval~\cite{Du} and Duval and Sibony~\cite{DS}, we
may construct non-orientable compact smooth submanifolds $M$ of $\cc^2$
which are rationally convex and the first homotopy groups of their
complements $\cc^2\setminus{M}$ are not trivial. Finally, a direct
application of algebraic topology yields the following extension of
Theorem~\ref{main}.  We present more applications in the third section
of this paper.

\begin{corollary}\label{homology}
Let $K$ be a compact rationally convex set in $\cc^n$, 
its complement is $(n-2)$-connected for $n\geq3$:
\begin{equation}\label{eq1}
\pi_k(\cc^n\setminus{K})=0\quad\hbox{when}\quad1\leq{k}\leq{n}-2;
\end{equation} 
and for any commutative group $G$,
\begin{equation}\label{eq2}
H_k(\cc^n\setminus{K},G)=0\quad\hbox{when}\quad1\leq{k}\leq{n}-2.
\end{equation} 
\end{corollary}

\begin{proof} We know that $\cc^n\setminus{K}$ is arcwise and 
simply connected, because Theorem~\ref{main}. Hence, applying 
Hurewicz theorem~\cite[p.~398]{Sp}, we have that equation~(\ref{eq1}) 
can be directly deduced from equation~(\ref{eq2}).

On the other hand, because of their definition, every rationally 
convex set $K$ has a system of Stein open neighbourhoods (open rational
polyhedra) $\{U_\beta\}$ in $\cc^n$, and so $K=\bigcap_\beta{U}_\beta$.  
It is well known that each Stein open set $U_\beta\subset\cc^n$ has the
homotopy type of a $CW$-complex of real dimension less than or equal to
$n$, see \cite[p.~26]{Di} or \cite[p.~548]{Ha}. Whence, all cohomology
groups $H^k(U_\beta,G)$ vanish for any commutative group $G$ and 
$k\geq{n+1}$. Alexander (and \v{C}ech) cohomology groups $\bar{H}^k(K,G)$ 
vanish as well for $k\geq{n}+1$, because they are equal to the direct 
limits of $H^k(U_\beta,G)$ when $U_\beta$ runs over a system of open 
neighbourhoods of $K$, see for example \cite[p.~291]{Sp} or 
\cite[p.~348]{Be}. Finally, a direct application of Alexander duality 
theorem yields equation~(\ref{eq2}), as we wanted, see \cite[p.~296]{Sp} 
or \cite[p.~352]{Be}. Notice that singular $H_k(\cdot)$ and reduced 
$\tilde{H}_k(\cdot)$ homology groups coincide for $k\geq1$.  
\end{proof}

Recall that a compact rational polyhedron $\Pi$ in 
$\cc^n$ is a set defined by the finite intersection,
\begin{equation}\label{eq3}
\textstyle{\Pi\,:=\,\overline{B}\cap\bigcap_kQ_k^{-1}(D_k),}
\end{equation}
where $\overline{B}\subset\cc^n$ is any compact ball, $\{Q_k\}$ is 
a finite collection of holomorphic polynomials on $\cc^n$, and each 
$D_k\subset\cc$ is either the compact unit disk $\{|w|\leq1\}$ or the 
closed ring $\{|w|\geq1\}$. Obviously, every compact rational polyhedron 
is rationally convex; and by an open rational polyhedron we understand 
the topological interior of a compact one.

We prove Theorem~\ref{main} in the next section, and several 
applications are presented in the third section of this paper. 

\section{Proof of main Theorem~\ref{main}}

The proof is based in the following two classical results on holomorphic
polynomials. Firstly, we need Lefschetz result on hyperplane sections,
where it is stated that the one-point compactification of a smooth
algebraic hypersurface is arcwise and simply connected. A simplified 
proof and several generalisations can be found in Hamm~\cite{Ha} or
Kar\v{c}jauskas~\cite{Ka}.

\begin{theorem}[Lefschetz]\label{lefschetz}
Let $X$ be a non-singular (smooth) algebraic hypersurface in $\cc^n$, for 
$n\geq3$. There then exists a finite real number $R_0>0$ such that the 
relative pair $(X,X\setminus{B}_R)$ is $(n-2)$-connected for every open 
ball $B_R\subset\cc^n$ of radius $R\geq{R_0}$ and centre in the origin.
\end{theorem}

This version of Lefschetz theorem can be easily deduced from the remark
in~\cite[p.~552]{Ha}, we just need to observe that $X\cap{B_R}$ is Stein
and to fix $H$ equal to the complement of $B_R$. Besides, the natural
existence of $R_0>0$ in Theorem~\ref{lefschetz} is easily understood 
after comparing Theorem~6.9 in \cite[p.~26]{Di} with Theorem~2 or~3 in
\cite{Ha}.  On the other hand, we also need the following result where 
it is shown that any non-constant holomorphic polynomial defined on 
$\cc^n$ almost induces a locally trivial fibre bundle on $\cc^n$, see
Verdier~\cite{Ve}, Broughton~\cite{Bo} and H\'a Huy Vui~\cite{HL}.

\begin{theorem}\label{fibration}
Let $Q:\cc^n\to\cc$ be a non-constant holomorphic polynomial. Then, there 
exists a finite set $\Lambda_Q$ in $\cc$ such that the fibres of $Q$ induce 
a locally trivial fibre bundle of $\cc^n\setminus{Q}^{-1}(\Lambda_Q)$ with 
base on $\cc\setminus\Lambda_Q$.
\end{theorem}

We need to combine previous two theorems into the following 
version of  Lefschetz theorem for \textit{fat} open sets.

\begin{lemma}\label{compact} 
Let $Q:\cc^n\to\cc$ be a holomorphic polynomial, for $n\geq3$, 
and $\Omega$ be a bounded open set in $\cc$ whose compact closure
$\overline\Omega$ is diffeomorphic to the unit disk $\{|w|\leq1\}$ 
and does not intersect the finite set $\Lambda_Q$ defined in
Theorem~\ref{fibration}. Given any ball $B_R$ in $\cc^n$ of finite 
radius $R>0$ and centre in the origin, there exists an unbounded
$(n-2)$-connected set $\Sigma_R$ in $\cc^n$ such that 
$B_R\cap\Sigma_R$ is equal to the open set $Q^{-1}(\Omega)\cap{B_R}$.  
\end{lemma}

\begin{proof} The result is trivial when $Q$ is a constant polynomial, 
so we suppose from now on that $Q$ is a non-constant function. The 
compact set $\overline\Omega$ is contractible in $\cc\setminus\Lambda_Q$, 
for it is diffeomorphic to the unit disk, so the fibres of $Q$ induce 
a global trivial fibration of $Q^{-1}(\overline\Omega)$ with base on 
$\overline\Omega$; recall Theorem~\ref{fibration} and consider 
\cite[p.~140]{AGP} or \cite[p.~27]{Di}. That is, given $w_0$ in $\Omega$, 
there exists a smooth function $g$ defined from $Q^{-1}(\overline\Omega)$ 
onto the non-singular fibre $X:=Q^{-1}(w_0)$, such that the pair $[Q,g]$ 
is a fibration diffeomorphism (it preserves fibres) from 
$Q^{-1}(\overline\Omega)$ onto $\overline\Omega\times{X}$.

In particular, the restriction $g|_X$ is an auto-diffeomorphism; and 
we may even fix $g$ such that $g|_X$ is the identity function. Besides, 
the norm $\|g(z)\|$ converges to infinity if and only if $\|z\|$ goes 
to infinity for $z$ in $Q^{-1}(\overline\Omega)$, recall that
$\overline\Omega$ is compact. Given any finite radius $R>0$, 
choose a real number $\rho\gg{R}$ such that the compact set
$g(Q^{-1}(\overline\Omega)\cap\overline{B_R})$ is contained 
inside $X\cap{B}_\rho$, and so define the bounded open set,
\begin{equation}\label{ER}
E_\rho\;\,:=\;[Q,g]^{-1}\big(\Omega\times(X\cap{B}_\rho)\big)
\;\subset\;Q^{-1}(\Omega). 
\end{equation} 
Where $[Q,g]^{-1}$ is the inverse diffeomorphism defined from
$\overline\Omega\times{X}$ onto $Q^{-1}(\overline\Omega)$. Since
$X=Q^{-1}(w_0)$ and the intersection $Q^{-1}(\Omega)\cap{B_R}$ is
contained inside $E_R$, because of the way we chose $\rho\gg{R}$, 
we can deduce that the following pair of equalities holds,
\begin{equation}\label{ES}\begin{array}{rcl}
E_\rho\cap{B}_R&=&Q^{-1}(\Omega)\cap{B}_R,\\
X\cup{E_\rho}&=&[Q,g]^{-1}\big([\{w_0\}\times
X]\cup[\Omega\times(X\cap{B}_\rho)]\big).  
\end{array}\end{equation} 
The open set $\Omega$ is contractible, for it is diffeomorphic to $\cc$.
Therefore, the spaces $X\cup{E}_\rho$ and $X$ have the same homotopy type.
Actually, it is easy to deduce that $X$ is a strong deformation retract 
of $X\cup{E}_\rho$.

On the other hand, Theorem~\ref{lefschetz} implies that the pair
$(X,X\setminus{B}_\sigma)$ is $(n-2)$-connected for any radius
$\sigma\gg\rho$ large enough. We may choose $\sigma\gg\rho$ such that
$X=Q^{-1}(w_0)$ meets transversally the sphere $\partial{B}_\sigma$, 
see Theorem~6.9 in \cite[p.~26]{Di}. Hence, if $B^{\com}_\sigma$ is the
complement of $B_\sigma$ in $\cc^n$, we automatically have that the 
pair $(B^{\com}_\sigma\cup{X},B^{\com}_\sigma)$ is $(n-2)$-connected.  
Finally, we may even fix $\sigma\gg\rho$ larger enough such that the 
open ball $B_\sigma$ contains the compact set $\overline{E_\rho}$ defined 
according to~(\ref{ER}). Thus, the spaces $X\cup{E}_\rho$ and $X$ have 
the same homotopy type and are equal outside the ball $B_\sigma$. The 
pair $(\Sigma_R,B^{\com}_\sigma)$ is then $(n-2)$-connected for 
$\Sigma_R$ defined as the union $B^{\com}_\sigma\cup{X}\cup{E}_\rho$.  
A direct application of the long exact sequence, see \cite[p.~87]{AGP} 
or \cite[p.~374]{Sp}
$$\to\pi_k(B^{\com}_\sigma)\to\pi_k(\Sigma_R)\to\pi_k(\Sigma_R,B^{\com}_\sigma)\to$$
automatically yields that $\Sigma_R$ is $(n-2)$-connected as well, recall
that the complement $B^{\com}_\sigma$ is $(2n-2)$-connected. It is easy 
to deduce that $\Sigma_R$ is unbounded and arcwise connected, because all 
the arcwise connected components of $X$ and $X\cup{E}_\rho$ are unbounded, 
and so they all meet the complement $B^{\com}_\sigma$. Finally, since
$\sigma\gg\rho\gg{R}$ and using equation~(\ref{ES}), we may conclude that 
$B_R\cap\Sigma_R$ is equal to the open set $Q^{-1}(\Omega)\cap{B}_R$, as 
we wanted.  
\end{proof}

The proof of Theorem~\ref{main} follows now an inductive process. 
Recall the definition of a compact rationally polyhedron given in
equation~(\ref{eq3}). The complement of any compact ball $\overline{B}$ in
$\cc^n$ is obviously simply connected. Thus, given any rationally convex
set $K$ in $\cc^n$ whose complement is simply connected, we only need to
prove that the union of $Q^{-1}(\Omega)$ and $\cc^n\setminus{K}$ is also
simply connected for all open sets $\Omega$ in $\cc$ and holomorphic
polynomials $Q$ in $\cc^n$. We prove previous statement in several steps.

\begin{lemma}\label{contractible} 
Let $Q:\cc^n\to\cc$ be a holomorphic polynomial, for $n\geq3$, and $\Pi$
be any compact rational polyhedron in $\cc^n$ whose complement is simply
connected. The union of $Q^{-1}(\Omega)$ and $\cc^n\setminus\Pi$ is then
simply connected for every bounded open set $\Omega$ in $\cc$ whose
compact closure $\overline{\Omega}$ is diffeomorphic to the unit disk
$\{|w|\leq1\}$ and does not intersect the finite set $\Lambda_Q$ defined
in Theorem~\ref{fibration}.  
\end{lemma}

\begin{proof} The result is trivial whenever $Q$ is a constant polynomial,
so we suppose from now on that $Q$ is a non-constant function. Define
$\Pi^{\com}$ to be the complement $\cc^n\setminus\Pi$, so the open set 
$\Pi^{\com}$ is simply connected because of the given hypotheses, and 
arcwise connected because $\Pi$ is rationally convex. Given an arbitrary 
closed path (loop) $\Upsilon$ inside $Q^{-1}(\Omega)$ union 
$\Pi^{\com}$, choose a ball $B_R$ with centre in the origin and 
radius $R>0$ large enough such that the compact sets $\Pi$ 
and $\Upsilon$ are both contained inside $B_R$.

Lemma~\ref{compact} implies the existence of an unbounded 
$(n-2)$-connected space $\Sigma_R$ in $\cc^n$ such that
$B_R\cap\Sigma_R$ is equal to the open set $Q^{-1}(\Omega)\cap{B}_R$.  
Notice in particular that $\Sigma_R$ is both arcwise and simply 
connected because $n\geq3$. Since $\Pi$ is contained in $B_R$, 
the union of $Q^{-1}(\Omega)$ and $\Pi^{\com}$ is equal to both 
spaces: $\Sigma_R\cup\Pi^{\com}$ and the complement of 
$\Pi\setminus{Q}^{-1}(\Omega)$ in $\cc^n$. We show that the loop
$\Upsilon$ is homotopically trivial in $\Sigma_R\cup\Pi^{\com}$ via 
an alternative version of Van-Kampen's theorem \cite[p.~63]{AGP}.

The spaces $\Sigma_R$ and $\Pi^{\com}$ are both simply connected 
because of Lemma~\ref{compact} and the given hypotheses, so we suppose 
that $\Upsilon$ meets $\Sigma_R$ and $\Pi^{\com}$. Since the loop 
$\Upsilon$ is contained in $B_R$, we can decompose it as the finite 
sum $\sum_k\Upsilon_k$ of non-necessarily closed paths $\Upsilon_k$ 
contained inside either the open set $B_R\cap\Sigma_R$ or the open 
complement $\Pi^{\com}$, and moreover, such that the end points of the 
paths $\Upsilon_k$ are all in the intersection $\Sigma_R\cap\Pi^{\com}$.  
Actually, since $\Pi$ is compact and $\Sigma_R$ is unbounded, the space 
$\Sigma_R\cap\Pi^{\com}$ is not empty, and so we can fix a point $w_0$ 
there.

Suppose that $\Sigma_R\cap\Pi^{\com}$ is also arcwise connected. We 
can then construct closed paths (loops) $T_k$ by joining both end points 
of each $\Upsilon_k$ to the fix point $w_0$ with arcs of the appropriated 
orientation in $\Sigma_R\cap\Pi^{\com}$. Every loop $T_k$ is contained 
inside one of the simply connected spaces $\Sigma_R$ or $\Pi^{\com}$. 
Hence, the original loop $\Upsilon$ is homotopic to the finite sum 
$\sum_kT_k$ in the union $\Sigma_R\cup\Pi^{\com}$, and moreover, the 
loops $T_k$ are all homotopically trivial there. In other words, the union 
of $Q^{-1}(\Omega)$ and $\Pi^{\com}$ is simply connected, because it 
is equal to $\Sigma_R\cup\Pi^{\com}$ and the original loop $\Upsilon$ is 
homotopically trivial there. We just need to show now that the intersection 
of $\Sigma_R$ and $\Pi^{\com}$ is indeed arcwise connected, in order to 
complete the proof of Lemma~\ref{contractible}.

Recall the proof of Corollary~\ref{homology}: Every rationally convex 
set $K$ in $\cc^n$ has a system of Stein open neighbourhoods, and so the
Alexander cohomology group $\bar{H}^{2n-2}(K)$ vanishes for $n\geq3$, 
see \cite[p.~26]{Di} and \cite[p.~291]{Sp}. The reduced homology group
$\tilde{H}_1(\cc^n\setminus{K})$ vanishes as well because Alexander
duality theorem \cite[p.~296]{Sp}. Since $\Sigma_R\cup\Pi^{\com}$ 
is the complement in $\cc^n$ of the rationally convex set
$\Pi\setminus{Q}^{-1}(\Omega)$, as it is indicated the second paragraph of
this proof, we have that $\tilde{H}_1(\Sigma_R\cup\Pi^{\com})$ vanishes.
On the other hand, the couple $\{\Sigma_R,\Pi^{\com}\}$ is excisive
because $\Sigma_R\cup\Pi^{\com}$ is equal to the union of the open sets
$\Pi^{\com}$ and $\Sigma_R\cap{B}_R$. Recall Lemma~\ref{compact} and 
that $B_R$ contains the compact set $\Pi$. Moreover, the intersection of
$\Pi^{\com}$ and $\Sigma_R$ is not empty because $\Sigma_R$ is unbounded.
Hence, the following exact Mayer-Vietoris sequence for reduced homology
holds, see \cite[p.~189]{Sp} or \cite[p.~229]{Be},
$$0=\tilde{H}_1(\Sigma_R\cup\Pi^{\com})\to\tilde{H}_0(\Sigma_R\cap\Pi^{\com})  
\to\tilde{H}_0(\Sigma_R)\oplus\tilde{H}_0(\Pi^{\com})\to.$$

Notice that a topological space $Y$ is arcwise connected if and only if
the reduced homology group $\tilde{H}_0(Y)$ vanishes. The fact that $\Pi$
is rationally convex and Lemma~\ref{compact} automatically imply that
$\Pi^{\com}$ and $\Sigma_R$ are both arcwise connected spaces, and so
their intersection $\Sigma_R\cap\Pi^{\com}$ is also arcwise connected. We 
can conclude that $Q^{-1}(\Omega)$ union $\Pi^{\com}$ is simply connected, 
after recalling a paragraph above. 
\end{proof}

The next steep implies extending Lemma~\ref{contractible} to consider
arbitrary open sets $\Omega$ in $\cc$ with a finite number of holes.

\begin{lemma}\label{holes}
Let $Q:\cc^n\to\cc$ be a holomorphic polynomial, for $n\geq3$, and $\Pi$
be any compact rationally polyhedron in $\cc^n$ whose complement is simply
connected. The union of $Q^{-1}(\Omega)$ and $\cc^n\setminus\Pi$ is then
simply connected for every open connected set $\Omega$ in $\cc$ with a
finite number of holes. 
\end{lemma}

The fact that $\Omega\subset\cc$ has a finite number of holes means 
that his fundamental group $\pi_1(\Omega)$ is free and finitely generated.  
Moreover, the conclusions of Lemma~\ref{holes} are quite interesting when
$\Pi$ is a compact ball.

\begin{proof} The result is trivial whenever $Q$ is a constant polynomial,
so we suppose from now on that $Q$ is a non-constant function. Define
$\Pi^{\com}$ to be the complement $\cc^n\setminus\Pi$, so the open set 
$\Pi^{\com}$ is simply connected because of the given hypotheses, and 
arcwise connected because $\Pi$ is rationally convex. Recall the finite set 
$\Lambda_Q$ defined in Theorem~\ref{fibration}. We prove this lemma by 
considering different open sets $\Omega_k$ in $\cc$ of increasing complexity.

Suppose that $\Omega_1$ is an open set in $\cc\setminus\Lambda_Q$ 
diffeomorphic to $\cc$.  Given any compact closed path (loop) $\Upsilon_1$ 
inside $Q^{-1}(\Omega_1)$ union $\Pi^{\com}$, the compact sets 
$Q(\Upsilon_1\cap\Pi)$ and $\Lambda_Q$ are \textit{far away} 
one from each other. Therefore, there exists a bounded open set 
$\Omega_2\subset\Omega_1$ such that the loop $\Upsilon_1$ is contained 
in $Q^{-1}(\Omega_2)$ union $\Pi^{\com}$ and the compact closure 
$\overline{\Omega_2}$ is both diffeomorphic to the unit disk $\{|w|\leq1\}$ 
and contained in $\Omega_1$. Notice that $\Lambda_Q$ does not meet 
$\overline{\Omega_2}$. A direct application of Lemma~\ref{contractible} 
yields that $\Upsilon_1$ is homotopically trivial in $Q^{-1}(\Omega_2)$ 
union $\Pi^{\com}$ and in the larger space $Q^{-1}(\Omega_1)$ union 
$\Pi^{\com}$. That is, the union of $Q^{-1}(\Omega_1)$ and $\Pi^{\com}$ 
is simply connected.

Suppose that $\Omega_3$ is any open connected set in 
$\cc\setminus\Lambda_Q$ with a finite number of holes. We may decompose 
$\Omega_3=V_1\cup{V}_2$ as the union of two open sets $V_k$ diffeomorphic 
to $\cc$ and whose intersection $V_1\cap{V}_2$ has a finite number of 
connected components. Every open set $Q^{-1}(V_k)\cup\Pi^{\com}$ is 
arcwise connected because his complement $\Pi\setminus{Q}^{-1}(V_k)$ 
is a rationally convex set. Previous paragraph implies that spaces 
$Q^{-1}(V_k)\cup\Pi^{\com}$ are both simply connected, for $k=1,2$. Finally, 
notice that $\bigcap_kQ^{-1}(V_k)$ union $\Pi^{\com}$ is arcwise connected 
because his complement is also a rationally convex set. We just need to apply 
Van-Kampen's theorem, in order to deduce that $Q^{-1}(\Omega_3)$ union 
$\Pi^{\com}$ is simply connected as well. That is, given two arcwise and 
simply connected open sets, their union is simply connected whenever their 
intersection is arcwise connected, see \cite[p.~63]{AGP} or \cite[p.~161]{Be}.

Finally, suppose that $\Omega_4$ is any open connected set in $\cc^n$ 
with a finite number of holes, so that $\Lambda_Q$ may intersect 
$\Omega_4$. Let $\Upsilon_4$ be any closed path (loop) in the open 
set $Q^{-1}(\Omega_4)$ union $\Pi^{\com}$.  There is a smooth loop 
$\Upsilon_5$ homotopic to $\Upsilon_4$ inside $Q^{-1}(\Omega_4)$ 
union $\Pi^{\com}$ such that $\Upsilon_5$ meets transversally each 
Whitney's strata of every singular fibre $Q^{-1}(s)$, for $s$ in 
$\Omega_4\cap\Lambda_Q$, recall Corollary~1.12 in \cite[p.~6]{Di} and 
section II-15 in \cite[pp.~114-118]{Be}. The loop $\Upsilon_5$ has real
dimension one and every fibre $Q^{-1}(s)$ has real codimension two in 
$\cc^n$, so the loop $\Upsilon_5$ meets no singular fibre $Q^{-1}(s)$ 
for $s$ in $\Omega_4\cap\Lambda_Q$. Previous paragraph implies 
that $\Upsilon_5$ is homotopically trivial in 
$Q^{-1}(\Omega_4\setminus\Lambda_Q)$ union $\Pi^{\com}$, 
because $\Omega_4\setminus\Lambda_Q$ has a finite number of 
holes. Thus, the original path $\Upsilon_4$ is homotopically trivial 
in the larger space $Q^{-1}(\Omega_4)$ union $\Pi^{\com}$, and 
so this union is simply connected, as we wanted to show. 
\end{proof}

Notice that previous lemma may be extended to consider any open set
$\Omega$ in $\cc$. Nevertheless, the proof becomes more complicated 
than what we really need in this paper. We are now in position to prove
Theorem~\ref{main}.

\begin{proof}\textbf{(Theorem~\ref{main})}. We begin by proving that 
the complement of a rational polyhedron is simply connected. Let
$\Pi\subset\cc^n$ be a compact rational polyhedron defined according 
to equation~(\ref{eq3}), for $n\geq3$. Notice that the complement of 
a compact ball $\cc^n\setminus\bar{B}$ is simply connected, because 
it has the homotopy type of the sphere $\mathcal{S}^{2n-1}$. From
Lemma~\ref{holes}, we have that the complement of the rational polyhedron
$Q_1^{-1}(D_1)\cap\bar{B}$ is also simply connected, because it is equal
to $Q_1^{-1}(\cc\setminus{D_1})$ union $\cc^n\setminus\bar{B}$. Following
an inductive process on the polynomials $Q_k$, and using Lemma~\ref{holes}
in every step, we can conclude that the complement $\cc^n\setminus\Pi$ is
simply connected as well.

On the other hand, let $K$ be any compact rationally convex 
set in $\cc^n$, and $\Upsilon$ a closed path in the complement
$\cc^n\setminus{K}$. We automatically have, from the definition, that
there exists of a compact rational polyhedron $\Pi\subset\cc^n$ which
contains $K$ and does not intersect $\Upsilon$. Thus, the path $\Upsilon$
is homotopically trivial in $\cc^n\setminus\Pi$, according to the previous
paragraph; and so, $\Upsilon$ is also homotopically trivial in the larger
set $\cc\setminus{K}$.  The complement of $K$ is then simply connected, 
as we wanted. 
\end{proof}

\section{Applications}

We want to finish this paper with the following application of
Corollary~\ref{homology}. Recall that the rational convex hull $r(K)$ of 
a compact set $K\subset\cc^n$ is equal to the intersection of all compact 
rational polyhedra $\Pi\subset\cc^n$ which contains $K$, see for example 
\cite{St} or \cite{AW}. The compact set $K$ is obviously contained in 
its rationally convex hull $r(K)$. Besides, given a continuous function 
$f$ defined from the $q$-dimensional sphere $\mathcal{S}^q$ into an open
subset $\Omega\subset\cc^n$, we say that $f$ is homotopically trivial 
in $\Omega$ if and only if it has a continuous extension to the compact
$(q+1)$-dimensional ball.

\begin{corollary} 
Let $K\subset\cc^n$ be a compact rationally convex set, for $n\geq3$, and 
$f$ be a continuous function defined from the sphere $\mathcal{S}^q$ into 
$\cc^n\setminus{K}$, which is not homotopically trivial there. The image 
$f(\mathcal{S}^q)$ intersects the rationally convex hull $r(K)$ whenever 
$1\leq{q}\leq{n-2}$. 
\end{corollary}

\begin{proof} If the image $f(\mathcal{S}^q)$ does not intersects $r(K)$, 
then, $f$ is homotopically trivial in the complement of $r(K)$ because 
Corollary~\ref{homology}. Hence, the function $f$ is also homotopically 
trivial in the larger set $\cc^n\setminus{K}$, contradiction. \end{proof}

We can deduce a similar result for homology. Recall that a $q$-cycle 
$C_q$ in an open subset $\Omega\subset\cc^n$ is formally a finite sum
$\sum_kg_kf_k$ of continuous functions $f_k$ defined from the standard
$q$-dimensional simplex $\Delta^q$ into $\Omega$, and coefficients 
$g_k$ in a commutative group $G$, see \cite{Sp} or \cite{Be}. We 
define the image of $C_q$ as the union of all $f_k(\Delta^q)$ for 
the coefficients $g_k\neq0$.

\begin{corollary} 
Let $K\subset\cc^n$ be a compact rationally convex set, for $n\geq3$, and 
$C_q$ be a $q$-dimensional cycle in the complement $\cc^n\setminus{K}$, 
which is not homologous to zero there. The image of $C_q$ intersects the 
rationally convex hull $r(K)$ whenever $1\leq{q}\leq{n-2}$. 
\end{corollary}

There exists a basic difference between homotopy and homology. Given a
continuous function $f$ defined from $\mathcal{S}^q$ into an open set
$\Omega\subset\cc^n$. We can express $\mathcal{S}^q$ as the union of two
compact hemispheres $D$ and $E$, and so we have that $C_q:=f|_D+f|_E$ is 
a $q$-cycle. The cycle $C_q$ is homologous to zero in $\Omega$, whenever 
$f$ is homotopically trivial there. Nevertheless, the cycle $C_q$ may be
homologous to zero, even when $f$ is not homotopically trivial.

Moreover, the difference between the first homotopy group $H_1(\cdot)$ and
the fundamental group $\pi_1(\cdot)$ is essential, even when $H_1(\cdot)$
is the abelianisation of $\pi_1(\cdot)$. For example, recalling that an
arc $\Upsilon$ is homeomorphic to the compact unit interval $[0,1]$ in the
real line, we automatically have that the Alexander cohomology groups
$\bar{H}^k(\Upsilon)$ vanish for every $k\geq1$, and so the first homology
group $H_1(\cc^n\setminus\Upsilon)$ vanishes for $n\geq1$ as well because
Alexander's duality theorem \cite[p.~352]{Be}. Nevertheless, Rushing
\cite[\S2.6]{Ru} produces several examples of arcs $\Upsilon$ in $\cc^n$
whose complement is not simply connected.  In particular, Rushing work
produces examples of arcs (and copies of the Cantor set) which cannot 
be rationally convex sets, according to Theorem~\ref{main}.

\bibliographystyle{plain}

\end{document}